\documentclass{amsart}

\newtheorem{theorem}{Theorem}[section]

\newtheorem{proposition}[theorem]{Proposition}

\theoremstyle{definition}

\theoremstyle{remark}

\numberwithin{equation}{section}
\renewcommand{\dim}{\mathrm{dim}}

\renewcommand{\span}{\mathrm{span}}

\newcommand{\conv}{\mathrm{conv}}

\newcommand{\N}{\mathbb{N}}

\newcommand{\X}{\mathbf{X}}

\newcommand{\Y}{\mathbf{Y}}
\newcommand{\Z}{\mathbf{Z}}
\newcommand{\B}{\mathbf{B}}

\renewcommand{\S}{\mathbf{S}}

\begin{document}

\title{On weakly extremal structures in Banach spaces}

\author{Jarno Talponen}
\address{University of Helsinki, Department of Mathematics and Statistics, Box 68, (Gustaf H\"{a}llstr\"{o}minkatu 2b) FI-00014 University 
of Helsinki, Finland}
\email{talponen@cc.helsinki.fi}

\date{\today}

\begin{abstract}
This paper deals with the interplay of the geometry of the norm and the weak topology in Banach spaces. 
Both dual and intrinsic connections between weak forms of rotundity and smoothness are discussed. Weakly locally uniformly rotund spaces, 
$\omega$-exposed points, smoothness, duality and the interplay of all the above are studied.   
\end{abstract}

\maketitle

\section{Introduction}
In introducing weaker forms of strong convexity in Banach spaces it is natural to apply the weak topology instead 
of the norm topology to describe the size of the extremal sections of the closed unit ball. In this article we study 
the weak geometry of the norm, that is, the interplay between the weak topology and the geometry of the norm.
While we investigate the duality between the smoothness and the convexity, it turns out
that also other types of interplay occur (see Theorem \ref{extreme}). 
Our study involves $\omega$-locally uniformly rotundity, $\omega$-exposedness and the G-smoothness of the
points $x\in\S_{\X}$, as well as the properties of the duality mapping $J\colon \S_{\X}\longrightarrow \mathcal{P}(\S_{\X^{\ast}})$. 
The concept of '$\omega$-exposedness' is natural in the sense that in a reflexive space $\X$ a point $x\in\S_{\X}$ is $\omega$-exposed 
if and only if it is exposed.

Let us mention the main themes appearing in this article. First, we give a characterization for 
$\omega$-exposed points $x\in\S_{\X}$ in terms of duality between the norm-attaining G-smooth functionals of the
dual space. This also leads to a sufficient condition for a point $x\in\S_{\X}$ to possess the $\omega$-LUR property. 
We also give a characterization for the reflexivity of a Banach space $\X$ in terms of an equivalent bidual $\omega$-LUR renorming 
of the bidual $\X^{\ast\ast}$. See e.g. \cite{Sullivan,90} for related research.

Secondly, we discuss the abundance of ($\omega$- or strongly) exposed points in various situations. For strongly exposed points this is 
a classical topic in the literature (see \cite{LP}) and well-studied especially in connection with the Radon-Nikodym property 
(see e.g. \cite[p.35]{JL}). In treating the question whether the $\omega$-exposed points are dense in $\S_{\X}$ 
we apply the topological properties of the duality mapping. 
 
We will use the following notations. Real Banach spaces are denoted by $\X,\Y$ and $\Z$ unless otherwise stated. 
We denote by $\B_{\X}=\{x\in\X:\ ||x||\leq 1\}$ the closed unit ball and by $\S_{\X}=\{x\in\X:\ ||x||=1\}$ the unit sphere. 
In what follows $\tau$ is a locally convex topology on $\X$.

For a discussion of basic concepts and results concerning the geometry of the norm we refer to \cite{HHZ} and
to the first chapter of \cite{JL}. 
The duality mapping $J\colon \S_{\X}\rightarrow \mathcal{P}(\S_{\X^{\ast}})$ is defined by
$J(x)=\{x^{\ast}\in\S_{\X^{\ast}}|x^{\ast}(x)=1\}$. If $A\subset \S_{\X}$ then we denote $J(A)=\bigcup_{a\in A}J(a)$. 
Denote by $\mathrm{NA}(\S_{\X^{\ast}})=J(\S_{\X})$ the set of norm-attaining functionals of $\S_{\X^{\ast}}$. 
Let us recall the following well known results:
\begin{theorem}\label{densesmooth}
If $\X$ is a separable Banach space, then the set of G-smooth points of $\S_{\X}$ is a dense $G_{\delta}$-set.
\end{theorem}
\begin{theorem}(Bishop-Phelps)\label{BPT}
The set $\mathrm{NA}(\S_{\X^{\ast}})$ is dense in $\S_{\X^{\ast}}$.
\end{theorem}
Recall that $\X$ is called an Asplund space if for any separable $\Y\subset\X$ it holds that $\Y^{\ast}$ is separable.
An Asplund space $\X$ satisfies that the F-smooth points of $\S_{\X}$ are a dense $G_{\delta}$-set (see \cite[Ch.1]{JL}). 
If $\X$ satisfies this conclusion for G-smooth points respectively, then $\X$ is called weakly Asplund.
A point $x\in\S_{\X}$ is called \emph{very smooth} if $x\in \S_{\X}\subset \S_{\X^{\ast\ast}}$ is G-smooth
considered in $\X^{\ast\ast}$.

If $f\in \S_{\X^{\ast}}$ is such that $f^{-1}(1)\cap \S_{\X}=\{x\}$ then $f$ is said to \emph{expose} $x$.
We say that $x\in \S_{\X}$ is a \emph{$\omega$-exposed} point if there is $f\in \S_{\X^{\ast}}$ such that whenever 
$(x_{n})\subset \B_{\X}$ is a sequence with $\lim_{n\rightarrow \infty}f(x_{n})=1$ then
$x_{n}\stackrel{\omega}{\longrightarrow}x$ as $n\rightarrow\infty$. If the same conclusion holds
for norm convergence then $x$ is called a \emph{strongly exposed} point. In such cases above $f$ is called
a \emph{$\omega$-} (resp. \emph{strongly}) \emph{exposing functional} for $x$.
We denote the $\omega$- (resp. strongly) exposed points of $\B_{\X}$ by $\omega\mathrm{-exp}(\B_{\X})$ 
(resp. $||\cdot||\mathrm{-exp}(\B_{\X})$).
A point $x\in\S_{\X}$ is called $\tau$-strongly extreme, 
if for all sequences $(z_{n}),(y_{n})\subset \B_{\X}$ such that 
$\frac{z_{n}+y_{n}}{2}\stackrel{||\cdot||}{\longrightarrow} x$ as $n\rightarrow\infty$ it holds that
$z_{n}-y_{n}\stackrel{\tau}{\longrightarrow}0$ as $n\rightarrow\infty$.
When $\tau$ is a locally convex topology on $\X$ we say that $x\in\S_{\X}$ is a $\tau$-Locally Uniformly Rotund point, 
$\tau$-LUR point for short, if for all sequences $(x_{n})\subset\B_{\X}$ such that $\lim_{n\rightarrow\infty}||x+x_{n}||=2$
it holds that $x_{n}\stackrel{\tau}{\longrightarrow}x$ as $n\rightarrow\infty$. If each $x\in\S_{\X}$ is $\tau$-LUR
then $\X$ is said to be $\tau$-LUR.
If $T$ is a topological space then a subset 
$A\subset T$ is called \emph{comeager} provided that it contains a countable intersection of subsets open and dense in $T$.

\section{Weak topology and convexity}
The most essential concept in this article is $\omega$-exposed point $x\in\S_{\X}$, which by its definition is exposed by a
$\omega$-exposing functional $f\in\S_{\X^{\ast}}$. Let us begin by characterizing the $\omega$-exposing functionals.
\begin{theorem}\label{weak}
Let $\X$ be a Banach space and suppose that $x\in \S_{\X},\ f\in \S_{\X^{\ast}}$ are such that $f(x)=1$. 
Then the following conditions are equivalent:
\begin{enumerate}
\item[(i)]{$f$ $\omega$-exposes $x$.}
\item[(ii)]{For each closed convex set $C\subset \B_{\X}$ such that $\sup_{y\in C}f(y)=1$ it holds that $x\in C$.}
\item[(iii)]{$f$ is a G-smooth point in $\X^{\ast}$, i.e. there is a unique $\psi\in \S_{\X^{\ast\ast}}$ such that $\psi(f)=1$.}
\end{enumerate}
\end{theorem}

Before giving the proof we will make some remarks. The above result must be previously known. For example 
by applying \cite[Thm.1,Thm.3]{ZZ} one can deduce the equivalence (i)$\iff$(iii). However, we will give a more elementary proof below.  
If above $\X$ is reflexive and $f$ exposes $x$ then one can see by the weak compactness of $\B_{\X}$ 
that actually $f$ $\omega$-exposes $x$.

\begin{proof}[Proof of Theorem \ref{weak}]
The equivalence (i)$\Leftrightarrow$(ii) follows easily by applying Mazur's theorem to $\conv(\{x_{n}|n\in\N\})$,
where $(x_{n})\subset \B_{\X}$ is a sequence such that $f(x_{n})$ tends to $1$ as $n\rightarrow\infty$.
 
Direction (iii)$\implies$(ii): 
Towards this suppose that $C\subset \B_{\X}$ is a closed convex subset so that
$\sup_{y\in C}f(y)=1$. Fix a sequence $(y_{n})\subset C$ satisfying $f(y_{n})\rightarrow 1$ as $n\rightarrow \infty$.
Observe that $x\in \S_{\X}\subset \S_{\X^{\ast\ast}}$ is the unique norm-one functional supporting $f$, since $f$ is a G-smooth
point of $\X^{\ast}$. By applying the \u Smulyan lemma to $x$ and $(y_{n})$ considered in $\X^{\ast\ast}$
we get that $y_{n}\stackrel{\omega^{\ast}}{\longrightarrow}x$ in $\X^{\ast\ast}$ as $n\rightarrow\infty$. This means that
$h(y_{n})\rightarrow h(x)$ as $n\rightarrow \infty$ for all $h\in \X^{\ast}$, so that $y_{n}\stackrel{\omega}{\longrightarrow}x$
in $\X$ as $n\rightarrow \infty$. Thus $x\in\overline{\conv}(\{y_{n}|n\in \N\})\subset \overline{\conv}^{\omega}(C)=C$
by Mazur's theorem, since $C$ is a norm-closed convex set. 

Conversely, suppose that (i) holds and let $\phi\in \B_{\X^{\ast \ast}}$ be an arbitrary point such that $\phi(f)=1$. 
We claim that $\phi=x$, which yields that $f$ is a G-smooth point. Let $g\in \X^{\ast}$ be arbitrary and recall that 
$\B_{\X}$ is $\omega^{\ast}$-dense in $\B_{\X^{\ast \ast}}$ by Goldstine's theorem. In particular, 
$\phi\in\overline{\B_{\X}}^{\omega^{\ast}}$. Thus 
$$\B_{\X}\cap \{\psi\in \B_{\X^{\ast\ast}}:\ |\psi(f)-1|<\frac{1}{n}\ \mathrm{and}\ |\psi(g)-\phi(g)|<\frac{1}{n}\}\neq\emptyset$$
for all $n\in \N$. Hence we may pick a sequence $(y_{n})\subset \B_{\X}$ for which $f(y_{n})\rightarrow 1$ and 
$g(y_{n})\rightarrow \phi(g)$ as $n\rightarrow \infty$. Since $f$ $\omega$-exposes $x$ we know that 
$y_{n}\stackrel{\omega}{\longrightarrow}x$ in $\X$ as $n\rightarrow\infty$. 
This yields that $\phi(g)=\lim_{n\rightarrow\infty}g(y_{n})=g(x)$ and that $\phi=x$ as this equality holds for all $g\in \X^{\ast}$.
\end{proof}

\begin{proposition}
Let $\X$ be a Banach space, $x\in\S_{\X}$ a F-smooth point and $f\in\S_{\X^{\ast}}$ a G-smooth point such that
$f(x)=1$. Then $x$ is a $\omega$-LUR point.
\end{proposition}
\begin{proof}
Suppose $(x_{n})\subset\B_{\X}$ is a sequence such that $||x_{n}+x||\rightarrow 2$ as $n\rightarrow \infty$.
By Theorem \ref{weak} the functional $f$ $\omega$-exposes $x$. Thus it suffices to show that
$f(x_{n})\rightarrow 1$ as $n\rightarrow \infty$.

By the Hahn-Banach Theorem one can find a sequence of functionals $(g_{n})\subset \S_{\X^{\ast}}$ such that
$g_{n}\left(\frac{x+x_{n}}{2}\right)=\left|\left|\frac{x+x_{n}}{2}\right|\right|$ for each $n\in\N$.
Clearly $g_{n}(x)\rightarrow 1$ and $g_{n}(x_{n})\rightarrow 1$ as $n\rightarrow \infty$.
Hence the F-smoothness of $x$ together with the \u Smulyan Lemma yields that 
$g_{n}\stackrel{||\cdot||}{\longrightarrow} f$ as $n\rightarrow\infty$. Thus we obtain that $f(x_{n})\rightarrow 1$ 
as $n\rightarrow\infty$.
\end{proof}  
\begin{proposition}
Let $\X$ be a Banach space and suppose that $x^{\ast}\in \S_{\X^{\ast}}$ is a very smooth point. Then there exists 
$x\in \S_{\X}\subset \S_{\X^{\ast\ast}}$ such that $x^{\ast}\in\S_{\X^{\ast}}\subset \S_{\X^{\ast\ast\ast}}$ $\omega$-exposes
$x$ in $\X^{\ast\ast}$.
\end{proposition}
\begin{proof}
Since by the definition $x^{\ast}$ is G-smooth in $\X^{\ast\ast\ast}$ it suffices to show that there exists $x\in \S_{\X}$
such that $x^{\ast}(x)=1$. Indeed, once this is established we may apply Theorem \ref{weak} to obtain the claim.

Let $x^{\ast\ast}\in \S_{\X^{\ast\ast}}$ be such that $x^{\ast\ast}(x^{\ast})=1$. By Goldstein's theorem $\B_{\X}\subset\B_{\X^{\ast\ast}}$
is $\omega^{\ast}$-dense. Pick a sequence $(x_{n})\subset\S_{\X}$ such that $x^{\ast}(x_{n})\rightarrow 1$ as $n\rightarrow \infty$. 
Observe that according to Theorem \ref{weak} the functional $x^{\ast}$ considered in $\S_{\X^{\ast\ast\ast}}$ $\omega$-exposes 
$x^{\ast\ast}$ in $\X^{\ast\ast}$. Hence $x_{n}\stackrel{\omega}{\longrightarrow} x^{\ast\ast}$ in $\X^{\ast\ast}$ 
as $n\rightarrow\infty$. Since $\X\subset \X^{\ast\ast}$ is $\omega$-closed by Mazur's theorem, we obtain that 
$x^{\ast\ast}\in \S_{\X}\subset\S_{\X^{\ast\ast}}$.  
\end{proof}

It is a natural idea to characterize reflexivity of Banach spaces in terms of suitable equivalent renormings
(see e.g. \cite{Hajek}). 
\begin{theorem}
The following conditions are equivalent:
\begin{enumerate}
\item[(1)]{$\X$ is reflexive.}
\item[(2)]{$\X$ admits an equivalent renorming such that $\X^{\ast\ast}$ is $\omega$-LUR.}
\item[(3)]{$\X$ admits an equivalent renorming such that $\Lambda\subset\S_{\X^{\ast\ast}}$ given by
\[\Lambda=\left\{\phi\in \S_{\X^{\ast\ast}}| \forall\ (\phi_{n})_{n\in\N}\subset\S_{\X^{\ast\ast}}:
\sup_{n}||\phi+\phi_{n}||=2\ \implies\ \phi\in[(\phi_{n})_{n\in\N}]\right\}\]
satisfies that $[\Lambda]=\X^{\ast\ast}$.}
\end{enumerate}
\end{theorem}
\begin{proof}
If $\X$ is reflexive then it is weakly compactly generated and hence admits an equivalent LUR norm, see e.g. \cite[p.1784]{JL2}.
Thus, by using reflexivity again we obtain that $\S_{\X^{\ast\ast}}$ is LUR.

Direction (2)$\implies$(3) follows by using Mazur's theorem that for convex sets weak and norm closure coincide. 
Since reflexivity is an isomorphic property we may assume without loss of generality in proving 
direction (3)$\implies$(1) that $\X$ already satisfies $[\Lambda]=\X^{\ast\ast}$.
Fix $\phi\in\Lambda$. Select a sequence $(f_{n})\subset \S_{\X^{\ast}}$ such that
$\phi(f_{n})\rightarrow 1$ as $n\rightarrow\infty$. Pick a sequence $(x_{n})_{n\in\N}\subset \S_{\X}$ such that 
$f_{n}(x_{n})\rightarrow 1$ as $n\rightarrow\infty$. This means that
\[||\phi+x_{n}||_{\X^{\ast\ast}}\geq (\phi+x_{n})(f_{n})\longrightarrow 2\ \mathrm{as}\ n\rightarrow\infty.\]
Hence by the definition of $\Lambda$ we obtain that $\phi\in [(x_{n})]$. Since $[\Lambda]=\X^{\ast\ast}$, this yields that 
$[\X]=\X^{\ast\ast}$ and hence $\X=\X^{\ast\ast}$ as $\X\subset\X^{\ast\ast}$ is a closed subspace.
\end{proof}

It turns out below that a smoothness property (namely Asplund) together with a weak convexity property 
(namely $\omega$-strongly extreme) yields in fact a stronger convexity property (namely the $\omega$-convergence), 
which is analogous to the '$\omega$-exposed situation'.

\begin{theorem}\label{extreme}
Let $\X$ be an Asplund Banach space and let $x\in\S_{\X}$ be a $\omega$-strongly extreme point. Suppose that 
$\{x_{n}|n\in\N\}\subset\B_{\X}$ is a set such that $x\in\overline{\conv}(\{x_{n}|n\in\N\})$. 
Then there is a sequence $(x_{n_{k}})_{k}\subset \{x_{n}|n\in\N\}$ such that $x_{n_{k}}\stackrel{\omega}{\longrightarrow}x$ as 
$k\rightarrow\infty$.
\end{theorem}
\begin{proof}
Consider convex combinations $y_{m}=\sum_{n\in J_{m}} a_{n}^{(m)}x_{n}\in\conv(\{x_{n}|n\in\N\})$
such that $y_{m}\rightarrow x$ as $m\rightarrow\infty$. Above $J_{m}\subset\N$ is finite and 
$a_{n}^{(m)}\geq 0$ are the corresponding convex weights for $m\in\N$. 

One shows easily that if $x_{n}^{\prime},x_{n}^{\prime\prime}\in \B_{\X},\ n\in\N,$ satisfy 
$\lambda_{n}x_{n}^{\prime}+(1-\lambda_{n})x_{n}^{\prime\prime}\rightarrow x$ with $\lambda_{n}\rightarrow \lambda\in (0,1]$
as $n\rightarrow\infty$, then $x_{n}^{\prime}\stackrel{\omega}{\longrightarrow} x$ as $n\rightarrow\infty$.

Fix $f\in \X^{\ast},\ \epsilon>0$ and put
\[K_{m}=\{n\in J_{m}|\ f(x_{n})<f(x)-\epsilon\}\quad \mathrm{for}\ m\in\N.\]
We claim that $\lambda_{m}=\sum_{n\in K_{m}}a_{n}^{(m)}\rightarrow 0$ as $m\rightarrow\infty$. Indeed, assume to the contrary that 
this is not the case. Then, by passing to a subsequence we may assume without loss of generality that 
$\lambda_{n}\rightarrow \lambda\in (0,1]$ as $n\rightarrow\infty$. Write $y_{m}=\lambda_{m}y_{m}^{\prime}+(1-\lambda_{m})y_{m}^{\prime\prime}$ for $m\in\N$, where
\[y_{m}^{\prime}=\frac{\sum_{n\in J_{m}}a_{n}^{(m)}x_{n}}{\lambda_{m}}\ \mathrm{and}\ y_{m}^{\prime\prime}=\frac{\sum_{n\in \N\setminus J_{m}}a_{n}^{(m)}x_{n}}{(1-\lambda_{m})}.\]
By the definition of the sequence $(y_{m}^{\prime})$ we have $f(y_{m}^{\prime})\leq f(x)-\epsilon$ for $m\in\N$ but this contradicts
the remark that $y_{m}^{\prime}\stackrel{\omega}{\longrightarrow} x$ as $m\rightarrow\infty$. 
 
Thus $\sum_{n\in K_{m}}a_{n}^{(m)}\rightarrow 0$ as $m\rightarrow\infty$ and a similar argument for 
$L_{m}=\{n\in J_{m}|f(x_{n})>f(x)+\epsilon\},\ m\in\N,$ gives that $\sum_{n\in L_{m}}a_{n}^{(m)}\rightarrow 0$ as $m\rightarrow\infty$.
These observations yield that 
\[\lim_{m\rightarrow\infty}\sum_{\substack{n\in J_{m}:\\ |f(x)-f(x_{n})|<\epsilon}}a_{n}^{(m)}=1\]
for any $\epsilon>0$. 

Since $f$ was arbitrary, we obtain that
\[\lim_{m\rightarrow\infty}\sum_{\substack{n\in J_{m}:\\ x_{n}\in U}}a_{n}^{(m)}=1,\]
where $U=\bigcap_{k}g_{k}^{-1}([g_{k}(x)-\epsilon,g_{k}(x)+\epsilon])$ and $g_{1},\ldots, g_{k}\in \X^{\ast},\ k\in\N$.
Recall that $x$ has a weak neighbourhood basis consisting of such sets $U$. In particular
\begin{equation}\label{eq: xom}
x\in\overline{\{x_{n}|n\in\N\}}^{\omega}.
\end{equation}

Let $\Y=\overline{\span}(\{x_{n}|n\in\N\})$. Since $\X$ is an Asplund space and $\Y$ is separable we obtain that $\Y^{\ast}$ is separable.
Fix a sequence $(h_{n})_{n}\subset \Y^{\ast}$, which is dense in $\Y^{\ast}$. Note that $(\B_{\Y},\omega)$ is metrizable
by the metric $d(x,y)=\sum_{n\in\N}2^{-n}|h_{n}(x-y)|,\ x,y\in\B_{\Y}$. Hence there is a sequence $(x_{n_{k}})_{k}\subset\{x_{n}|n\in\N\}$
such that $x_{n_{k}}\stackrel{\omega}{\longrightarrow} x$ as $k\rightarrow\infty$ in $\Y$. By the Hahn-Banach extension of $\Y^{\ast}$ 
to $\X^{\ast}$ it is straight-forward to see that $x_{n_{k}}\stackrel{\omega}{\longrightarrow}x$ as $k\rightarrow\infty$ also in $\X$. 
\end{proof}

\subsection{Density of $\omega$-exposed points}
\ \newline
Recall the following result due to Lindenstrauss and Phelps \cite[Cor.2.1.]{LP}: \\
\textit{If $C$ is convex body in an infinite dimensional separable reflexive Banach space, then the extreme points
of $C$ are not isolated in the norm topology.}\\
This result is the starting point for the studies in this section. 

Let us first consider a natural property of the duality mapping $J\colon \S_{\X}\rightarrow\mathcal{P}(\S_{\X^{\ast}})$. The following 
fact is an elementary topological statement about $J$ and it is proved here for convenience.
\begin{proposition}\label{equivprop}
The following conditions are equivalent for a Banach space $\X$:
\begin{enumerate}
\item[(1)]{For each relatively open non-empty $U\subset\S_{\X}$ the set 
$J(U)$ contains an interior point relative to $\mathrm{NA}(\S_{\X^{\ast}})$.}
\item[(2)]{For each relatively dense $A\subset \mathrm{NA}(\S_{\X^{\ast}})$ the subset 
$\{x\in\S_{\X}|\ J(x)\cap A\neq\emptyset\}\subset \S_{\X}$ is dense.}
\end{enumerate}
\end{proposition}
\begin{proof}
Suppose that (1) holds. Let $A\subset \mathrm{NA}(\S_{\X^{\ast}})$ be a dense subset. Assume to the contrary that 
$\{x\in\S_{\X}|\ J(x)\cap A\neq\emptyset\}\subset \S_{\X}$ is not dense. Thus there is
a non-empty set $U\subset \{x\in\S_{\X}|\ J(x)\cap A=\emptyset\}$, which is open in $\S_{\X}$. Hence
$J(U)\cap A=\emptyset$ where $J(U)\subset \S_{\X^{\ast}}$ is open according to condition (1). This contradicts the assumption that 
$A\subset \mathrm{NA}(\S_{\X^{\ast}})$ is dense and consequently condition (2) holds. 

Suppose that (2) holds and $U\subset \S_{\X}$ is a non-empty open set. Assume to the contrary that $J(U)$ does not contain an interior 
point relative to $\mathrm{NA}(\S_{\X^{\ast}})$. Then $\mathrm{NA}(\S_{\X^{\ast}})\setminus J(U)$ is dense in 
$\mathrm{NA}(\S_{\X^{\ast}})$. Hence condition (2) states that
$\S_{\X}\setminus U=J^{-1}(\mathrm{NA}(\S_{\X^{\ast}})\setminus J(U))$ is dense in $\S_{\X}$. This contradicts the 
assumption that $U$ is open and hence condition (1) holds.
\end{proof}

When $\X$ satisfies the above equivalent conditions of Proposition \ref{equivprop}, we say that $\X$ satisfies ($\ast$) for the sake of 
brevity. The following result describes this condition.

\begin{proposition}\label{*prop}
Suppose that $\X^{\ast}$ is an Asplund space. Then $\X$ satisfies condition $(\ast)$ if and only if 
the set of strongly exposed points of $\B_{\X}$ is dense in $\S_{\X}$.
\end{proposition} 
\begin{proof}
The 'if' case. Fix non-empty $A\subset \mathrm{NA}(\S_{\X^{\ast}})$. Suppose $\S_{\X}\setminus \{x\in\S_{\X}|\ J(x)\cap A\neq\emptyset\}$ 
contains an interior point relative to $\S_{\X}$. We aim to show that in such case $A\subset\S_{\X^{\ast}}$ is not dense.
Indeed, by the density assumption regarding the strongly exposed points we obtain that there is 
$x\in ||\cdot||\mathrm{-exp}(\B_{\X})$, which is not in the closure of $\{x\in\S_{\X}|\ J(x)\cap A\neq\emptyset\}$. 
If $f\in\S_{\X^{\ast}}$ is a strongly exposing functional for $x$ then there is $\epsilon>0$ such that 
$f(\{x\in\S_{\X}|\ J(x)\cap A\neq\emptyset\})\subset [-1,1-\epsilon]$. This gives that $||f-g||\geq \epsilon$ for all
$g\in A$, so that $A\subset\S_{\X^{\ast}}$ is not dense.

The 'only if' case. Let us apply the fact that $\X^{\ast}$ is Asplund. Denote by $\mathrm{F}$ the set of 
F-smooth points $x^{\ast}\in \S_{\X^{\ast}}$. Recall that $F\subset\S_{\X^{\ast}}$ is dense since 
$\X^{\ast}$ is Asplund. Condition ($\ast$) of $\X$ gives that $\{x\in\S_{\X}|\ J(x)\cap \mathrm{F}\neq\emptyset\}$ is dense in $\S_{\X}$.
By applying the \u Smulyan Lemma it is easy to see that each F-smooth functional $x^{\ast}\in \S_{\X^{\ast}}$ is norm-attaining and 
in fact a strongly exposing functional.
\end{proof}

The following main result is a version of the above-mentioned result by Lindenstrauss and Phelps.

\begin{theorem}\label{th:isolated}
Let $\X$ be a Banach space, which satisfies $\dim(\X)\geq 2$ and suppose that the following conditions hold:
\begin{enumerate}
\item[(1)]{$\mathrm{NA}(\S_{\X^{\ast}})\subset \S_{\X^{\ast}}$ is comeager.}
\item[(2)]{$\X^{\ast}$ is weakly Asplund.}
\end{enumerate} 
Then there does not exist a G-smooth point $x\in\S_{\X}$, which is $\omega$-isolated in $\omega\mathrm{-exp}(\B_{\X})$.\\ 
Moreover, if $\X$ satisfies additionally $(\ast)$, then $\omega\mathrm{-exp}(\B_{\X})\subset\S_{\X}$ is dense.
\end{theorem}
To comment on the assumptions shortly, observe that condition (1) above holds for instance if 
$\X$ has the RNP (see \cite{Phelps} and \cite[Thm.8]{Bour}) and $\X^{\ast}$ is weakly Asplund for instance if $\X^{\ast}$ is separable
by Theorem \ref{densesmooth}. It follows that the Asplund property of $\X^{\ast}$ is sufficient for both the conditions (1) and (2) 
to hold. Observe that we do not require $\X$ above to be separable, nor infinite dimensional. On the other hand, the assumption about 
the G-smoothness of $x$ can not be removed. For example consider $\ell^{\infty}_{n}$ for $n\in\{2,3,\ldots\}$.

\begin{proof}[Proof of Theorem \ref{th:isolated}]
According to the weak Asplund property of $\X^{\ast}$ there is a dense $G_{\delta}$-set of G-smooth points in $\S_{\X^{\ast}}$. 
We apply this fact together with assumption (1) as follows. By using the Baire category theorem we obtain that 
\[\mathcal{NG}=\{x^{\ast}\in \mathrm{NA}(\S_{\X^{\ast}})|x^{\ast}\ \mathrm{is\ G-smooth}\}\] 
is dense in $\S_{\X^{\ast}}$. Observe that by Theorem \ref{weak} all $f\in \mathcal{NG}$ are in fact $\omega$-exposing functionals.

Now, assume to the contrary that $x$ is a G-smooth $\omega$-isolated point in $\omega\mathrm{-exp}(\B_{\X})$. 
Then $x$ is $\omega$-exposed by a unique support functional $f\in \mathcal{NG}$. 
It is easy to see that since $f$ is a $\omega$-exposing and $\omega$-isolated functional, there is $\epsilon>0$ such that 
$f(\omega\mathrm{-exp}(\B_{\X})\setminus \{x\})\subset [-1,1-\epsilon]$.
Consequently, by the uniqueness of $f$ we obtain that $||f-g||\geq \epsilon$ for all $g\in \mathcal{NG}\setminus \{f\}$. 
Thus we obtain that the relatively open set
$\{h\in\S_{\X^{\ast}}:\ 0<||f-h||<\epsilon\}\subset\S_{\X^{\ast}}$ is non-empty as $\dim(\X)\geq 2$ and it does not intersect 
the dense subset $\mathcal{NG}\subset\S_{\X^{\ast}}$, which provides a contradiction. Hence the first part of the claim holds.

Finally, let us assume that $\X$ satisfies ($\ast$). Hence $\{x\in \S_{\X}|J(x)\cap \mathcal{NG}\neq\emptyset\}=J^{-1}(\mathcal{NG})$ 
is dense in $\S_{\X}$. Recall that the  points in $J^{-1}(\mathcal{NG})$ are $\omega$-exposed.
\end{proof}
 
\subsection*{Acknowledgments}
I thank the referee for many useful suggestions. This article is part of the writer's ongoing Ph.D. work, 
which is supervised by H.-O. Tylli. The work has been supported financially by the Academy of Finland projects 
\# 53968 and \# 12070 during the years 2003-2005 and by the Finnish Cultural Foundation in 2006.

\end{document}